\documentclass[11pt]{article}

\usepackage[ansinew]{inputenc}
\usepackage[psamsfonts]{amssymb}
\def\bH{I\!\!H} 
\def\bk{I\!\!k} 
\def\mm{{\mathfrak M}} 

\title{Frobenius Rational Loop Algebra}
\author{David Chataur and  Jean-Claude Thomas}

\begin{document}
\maketitle

    \begin{abstract} 
Recently R. Cohen and V. Godin have proved that the loop homology  $\bH_\ast 
(LM;\bk)$ of
a closed oriented manifold $M$ with coefficients in a field $\bk$  has the 
structure of a
unital Frobenius algebra without counit. In this paper  we  describe explicitly 
the
  dual of the  coproduct $ H_\ast (LM ;{\mathbb Q}   ) \to \left(H_\ast  (LM ; 
{\mathbb
Q}  )
\otimes H_\ast (LM; {\mathbb Q} )\right)_{\ast-m}
$  and  prove  it    respects
the Hodge decomposition.   
  
\end{abstract}

    \vspace{5mm}\noindent {\bf AMS Classification} : 55P35, 54N45,55N33,
  17A65,
    81T30, 17B55

    \vspace{2mm}\noindent {\bf Key words} : Frobenius algebra,  
topological quantum field theory, free loop space,  loop homology, rational 
homotopy,
Hochschild cohomology, Hodge decomposition.

\vspace{5mm}

 \centerline{\bf Introduction}

\vspace{3mm}

Let $\bk $ be a field with a unit element denoted 1. Singular homology 
(respectively 
cohomology) of a space $X$  with coefficients in $\bk$ is denoted $H_\ast(X)$ 
(respectively 
$H^\ast(X)$). For simplicity we identify  $H^\ast(X^{\times 2})$ with  $H^\ast(X)
^{\otimes 2} $ and $ H^\ast( X)$ with the graded dual of the homology 
$\left( H_\ast( X)\right)^\vee$.

    Let $M$ be a 1-connected closed  oriented $m$-manifold and let 
$M ^{S^1} $ be  the space of free loops. The space $M^{S^1}$      can be   
replaced
by a Hilbert manifold denoted $LM$ (see \cite{Ch} for more details) so that the 
evaluation map 
$p_0 :LM \to M\,, \quad \gamma \mapsto \gamma(0)$ is a smooth map.

 Chas and Sullivan \cite{C-S1} have constructed  a  natural product, called the 
{\it loop
product}
$$
P: H_\ast( LM \times LM) \to H_{\ast-m}(LM)
$$
so that $\bH_\ast(LM) := H_{\ast +m} (LM)$ is a commutative 
graded algebra. Later on   Chas and Sullivan, \cite{C-S2}, and Sullivan, 
\cite{S2}, have
introduced  ``higher string structures''.  Cohen-Godin, \cite{C-G},
have proved that some of  these higher string structures  can be defined in terms 
of two
dimensional positive boundary TQFT (topological quantum field theory). In, 
particular
they prove:

\vspace{3mm}
\noindent{\bf Theorem CG} (\cite{C-G}-Theorem 1) {\it There exists a loop 
coproduct
$$
\Phi : H_\ast (LM) \to H_{\ast -m}(LM\times LM)$$
such that $P$ and $\Phi$ define a structure of a  Frobenius algebra without 
counit on
$\bH_\ast(LM)$ .}

\vspace{3mm}

Recall that a  {\it   Frobenius algebra} is equivalent to a two dimensional 
positive
boundary TQFT, \cite{C-G}. In particular, it is  a unital associative and 
commutative
graded  algebra $A$ together  with a  (non zero) coassosiative coproduct   $\Phi 
: A \to
A
\otimes A$ which is a (degree d) $A$-linear map. Here the coproduct $\Phi$ is 
assumed to
be without   counit, compare with \cite{A1}, \cite{A2}. 
\vspace{3mm}

It results from the very definition of the loop coproduct and of the loop product 
that
the  Frobenius algebra structure on $H_\ast(LM)$ is natural with respect to 
 any orientation preserving smooth map   between   connected closed 
oriented
$m$-manifolds (see $\S3$-Proposition 1). Roughly speaking the above theorem 
asserts that
the functor
$ H_\ast( -)\circ L$ carries  the Poincar\'e duality structure on $H_\ast (M)$ to 
a
Frobenius algebra structure on $H_\ast(LM)$. Let us also notice that the 
orientation
preserving condition is necessary because the image of the fundamental class of 
$M$ in
$H_\ast(LM)$ is the unit of the loop product.

Our first result  is : 

\vspace{3mm}
\noindent{\bf Theorem A}   {\it If $\bk$ is of characteristic zero, the coproduct 
of
the   Frobenius algebra $H_\ast (LM) $ is  explicitly described    in terms  of 
Sullivan
models and with the Euler class of the   diagonal
 embedding $\Delta : M \to M \times M $. (See Theorem 10 in $\S 4$ for a more 
precise
statement.) 

Moreover, any orientation preserving smooth map $\varphi$  between  1-connected 
closed 
oriented
$m$-manifolds  such that $ H ^\ast(\varphi)$ is an isomorphism  induces an
isomorphism $ H _\ast(L\varphi)$ which respects the loop 
products and  the loop coproducts.}

 \vspace{3mm}
The loop product has already 
been  described in terms of Sullivan models in \cite{FTV2}-Theorem A.  For any 
field of
 coefficients    the loop product can be identified to the cup
product of the  Hochschild cohomology $H\! H^\ast (C^\ast M;C^\ast M)$ when 
$C^\ast M $
denotes the singular cocha\^\i ns on $M$,  \cite{C-J}   (\cite{Me}  when   $\bk
= {\mathbb R}$ and \cite{FTV2} when   $\bk
= {\mathbb Q}$). While the Hochschild cohomology of a differential
graded algebra,
$H\! H^\ast (A;A)$, is not natural with respect to graded algebra homomorphisms, 
we
know that if $ \varphi : A \to B$ is a quasi-isomophism of differential
graded algebras then there is  an isomorphism of graded algebras $H\!H(A;A)
\to H\! H(B;B)$, \cite{FMT}. Since the identification of $ H_\ast(LM)$ with $H\! 
H^\ast
(C^\ast M;C^\ast M)$ depends a priori of the Poincar\'e duality of $H_\ast(M)$  
it is
not possible to conclude from these results that  the abstract algebra structure 
on
$H_\ast(LM)$ is a homotopy invariant. (The difficulty here is to obtain a 
convenient 
chain level  representative of the  cap product with the orientation class, see
\cite{FTV1}-Theorem 13 and \cite{FTV2}).  After a preliminary version has been 
posted,
we  have learned   that  R. Cohen, J. Klein and D. Sullivan have announced at 
several
summer conferences in 2003 that the string topology operations are integral 
homotopy
invariant. (A preprint entitled ``On the homotopy invariance of string topology'' 
is in
preparation.)

\vspace{3mm}

 The tools  involved in the proof of theorem A  are 
Sullivan models \cite{S1}. In this context we naturally work with the
cohomology of  $LM$ instead of homology. Therefore we are interested by the {\it 
dual of
the  loop copoduct} $ \Phi^\vee  :
 H^\ast (LM\times LM)  \rightarrow  H^{\ast+m} (LM) \,,
$
(respectively  {\it the dual of the loop product} $ P^\vee  : H^\ast (LM) 
\rightarrow 
 H^{\ast+m} (LM\times LM)  \mbox{ )}$
rather
than the loop coproduct $\Phi$ (respectively  the loop product $P$). 

 Theorem A allows to do explicit examples. For instance we have  established in 
$\S$ 5
the  formula for
$\Phi ^\vee $ when $M={\mathbb C}P^n$ and we 
show  that the composition $
 \Phi \circ \mbox{loop product} :  (\bH(LM)\otimes\bH(LM))_\ast \to
(\bH(LM)\otimes\bH(LM))_\ast 
$ is not the trivial map.  

A second application of Theorem A concerns the ``
Hodge decomposition of
$H^\ast(LM)$''.

 Recall  that $LM$ is equipped with a natural $S^1$-action. Power maps $\varphi 
_N
:LM \to LM $ are then induced by the covering maps $ e^{i\theta} \mapsto 
e^{ni\theta}$.
The power operations $H^\ast (\varphi _n) : H^\ast (LM) \to H^\ast (LM)$ define, 
if $\bk$
is of characteristic zero,  a natural  Hodge decomposition, \cite{B-V} :
  $$
H^\ast (LM)= \oplus _{i\geq 0} H^\ast_{(i)}
\,, \mbox{ where } 
\left\{
\begin{array}{lll}
H^\ast _ {(0)} &= \tilde{H}^\ast(M)\\

H^n_{(i)}& = 0 \,, \mbox{ if } i>n\\

H^n_{(i)} &\mbox{ is an eingenspace of } H^n(\varphi_n)\\
& \mbox{ for the eigenvalue } n^i, 0<i\leq n\\
\end{array}
\right. \,.
$$

\vspace{5mm}
\noindent{\bf Theorem B.} {\it If $\bk$ is of characteristic zero the   dual  of 
the  
loop  coproduct $\Phi ^\vee $   respects the  Hodge
decomposition in $H^\ast(LM)$.}

\vspace{3mm}

 The paper is organized as follows;

1 - Definition of the  loop coproduct and of the loop product.

2 - Shriek map of a $k$-embedding.

3 - Sullivan representative of $\delta_{in}$ and of $\delta_{out}$.

4 -  Proof of theorem A.

5 -  Example $M= {\mathbb C}P^n$.

6 -  Proof of theorem B.

\vspace{3mm}
\noindent{\bf Acknowledgments.} We would like to thanks  R. Cohen, Y. Felix and 
J.
Stasheff for helpful comments.

\vspace{5mm}

 \centerline{\sc $\mathbf{\S 1}$ - Definition of the  loop coproduct and of the 
loop
product.}

 The purpose of this section is to  provide us with  an elementary definition
of the dual of the loop coproduct  and of the dual of the  loop product in 
algebraic topology.  

\vspace{3mm}

 Following \cite{S2} and \cite{C-G} the loop product (respectively  the loop
coproduct) is a  particular case of an ``operation''
$$
\mu_\Sigma : H_\ast(LM )^{\otimes p}  \to  H_\ast(LM )^{\otimes q}
$$
with $p=2\,,q=1$ (respectively  $p=1\,, q=2$) that will define a positive 
boundary TQFT. Here
$\Sigma$ denotes an oriented surface of genus 0 with a fixed parameterization of 
the 
$p+q$ boundary components:
$$
\left( S^1\right) ^{\coprod p}\stackrel{in}{\to} \partial _{in}\Sigma  \ 
\hookrightarrow
\Sigma
\ \hookleftarrow \partial _{out}\Sigma \stackrel{out}{\leftarrow} 
\left( S^1\right) ^{\coprod q}\,.
$$ 
Applying the  functor $Map(-,M)$ on gets the diagram
$$
(1) \qquad 
\left(LM\right)^{\times p} \stackrel{in ^\#}{\longleftarrow} Map(\Sigma, M) 
\stackrel{out ^\#}{\longrightarrow} \left(LM\right)^{\times q}
\,.
$$

Diagram (1) is homotopically equivalent to the diagram 

$$
(2) \qquad 
\left(LM\right)^{\times p} \stackrel{\rho_ {in} }{\longleftarrow} Map(c, M) 
\stackrel{\rho_{out}}{\longrightarrow} \left(LM\right)^{\times q}
\,.
$$
where $c$ denotes a reduced  Sullivan chord diagram with markings which is 
associated to
$\Sigma$
\cite{C-G}-$\S 2$.

\vspace{3mm}

In order to define the loop coproduct we restrict to the case $p=1$
and
$q=2$. Then 
$\Sigma $ is the oriented surface of genus 0 with one incoming  and two outcoming
components of the boundary. The associated Sullivan chord diagram with markings 
is
determined by the pushout diagram
$$
\begin{array}{ccc}
\ast \coprod \ast & \to &\ast\\
\downarrow&& \downarrow\\
S^1& \stackrel{in}{\to}& c
\end{array}
\,.
$$
In particular,  $c\simeq S^1\vee S^1$, $S^1 \stackrel{in}{\to} c$ is homotopic to 
the
folding map $\nabla : S^1 \to S^1\vee S ^1$ and $S^1\coprod S^1  
\stackrel{out}{\to} c$
is homotopic to the natural projection $S^1 \coprod S^1 \to S^1 \vee S^1$.

Diagram (2) then reduces to   a commutative  diagram of fibrations
$$
(3) \qquad 
\begin{array}{ccccccl}
 LM ^{\times 2} &\stackrel{\delta_{out}  }{\longleftarrow}&  LM \times_M 
LM&\stackrel{
\delta_{in} }{\longrightarrow}&  LM \\
\hspace{- 9mm} p_0 ^{\times 2}  \downarrow && \hspace{- 4mm}p_0 \downarrow &&
\hspace{2mm}\downarrow   (p_0, p_{\frac{1}{2}})
\\
M^{\times 2}  &\stackrel{\Delta  }{\leftarrow}&M
&\stackrel{\Delta}{\to}& M^{ \times 2}   
\end{array}
$$
where   

\noindent i) $p_t(\gamma) = \gamma(t)$,  $\gamma \in M^{(a,b)}$, $t\in 
(a,b)\subset 
I:=(0,1)$ 

\noindent ii) $\Delta (x) = (x,x)$ \\

\noindent iii) The fibration $LM \stackrel{(p_0, p_\frac{1}{2})} \to M^{\times 2} 
$ is
defined by the pullback diagram  
$$ (4) \qquad \qquad 
\begin{array}{lcl}
 LM &\stackrel{\delta' }{\to}& {M^{[0,\frac{1}{2}]}}\times {M^{[\frac{1}{2}, 
1]}}\\
\hspace{-9mm} (p_0 , p_{\frac{1}{2}})  \downarrow   &&\hspace{1mm} \downarrow
p_0\times p_\frac{1}{2}\times p_\frac{1}{2}\times p_1
\\
 M^{\times 2}  &\stackrel{\Delta'  }{\rightarrow}& M^{\times 4}
\end{array}\,, 
$$
and $\Delta' (x,y) = (x,y,y,x)$, $x,y \in M$

\noindent iv) the   squares are pullback diagrams of fibrations and  
$\delta_{in} $, $\delta _{out}$ are  the natural inclusion.\\

\vspace{3mm}

Diagram (3) converts into  the commutative diagram of $\bk$-linear maps
$$
\begin{array}{lccccccc}
H^\ast( LM ^{\times 2}) &\stackrel{H^\ast(
\delta _{out})}{\rightarrow}& H_\ast( LM\times _M LM) &\stackrel{{\delta_{in}}
^!}{\rightarrow}& H^\ast( LM)
\\
 H^\ast(p_0 ^{\times 2})
\uparrow &&
\hspace{- 9mm} H^\ast(p_{ 0}) \uparrow &&
\hspace{6mm}  
\uparrow H^\ast(p_0, p_{\frac{1}{2}}) \\

H^\ast( M^{\times 2})  &\stackrel{ H^\ast(\Delta) 
}{\rightarrow}& H^\ast(M) &\stackrel{\Delta ^!}{\rightarrow}& H^\ast(M^{ \times 
2})
\end{array}
$$
where $ {\delta_{in}}^!$ (respectively  $\Delta ^!$)  denotes the shriek map of 
the $m$-embedding
$\delta_{in}$ (respectively  of the diagonal embedding of $M$) 
(See section 2 for a precise definition).  By definition,  (\cite{C-G}- 
Definition 3),  the {\it
dual of the loop coproduct}  is 
$$
\Phi^\vee := \delta^!_{in} \circ H^\ast(\delta_{out}) :  \left(
H^\ast(LM)\right)^{\otimes 2} \to  
H^\ast(LM)\,.
$$

 To define the dual of the  loop product  one consider the diagram 
$$
(5) \qquad \begin{array}{ccccccl}
 LM ^{\times 2} &\stackrel{\delta_{out} }{\longleftarrow}&  LM \times_M
LM&\stackrel{Comp}{\longrightarrow}&  LM \\
\hspace{- 9mm} p_0 ^{\times 2}  \downarrow && \hspace{- 4mm}p_0 \downarrow &&
\hspace{2mm}\downarrow   p_0
\\
M^{\times 2}   &\stackrel{\Delta}{\longleftarrow} &M
&=& M    
\end{array}
$$
where $Comp$ denotes composition of free loops and $\delta_{out}$ is a fibrewize
embedding. Diagram (7) converts into  into  the commutative diagram of 
$\bk$-linear maps
$$
\begin{array}{lccccccc}
H^\ast( LM ^{\times 2}) &\stackrel{
\delta _{out}^!}{\leftarrow}& H_\ast( LM\times _M LM) 
&\stackrel{H^\ast(Comp)}{\leftarrow}&
H^\ast( LM)
\\
 H^\ast(p_0 ^{\times 2})
\uparrow &&
\hspace{- 9mm} H^\ast(p_{0}) \uparrow &&
\hspace{6mm}  
\uparrow H^\ast(p_0) \\

H^\ast( M^{\times 2})  &\stackrel{ \Delta^!  
}{\leftarrow}& H^\ast(M) &=& H^\ast(M)
\end{array}
$$
where $ {\delta_{out}}^!$   denotes the shriek map of the $m$-embedding
$\delta_{out}$  (See section 2 for a precise definition).  Following \cite{C-J} 
or \cite{C-G},
the  dual of the loop  product is defined by:
$$
P^\vee   = {\delta_{out} }^! \circ H^\ast (Comp)  : H^\ast(LM) \to  \left(
H^\ast(LM)\right)^{\otimes 2}\,.
$$

\vspace{5mm}

 \centerline{\sc $\mathbf{\S 2}$ - Shriek map of a $k$-embedding.}

\vspace{3mm}

{\it A fiber bundle pair}   
$$
(p,p') : (F,F') \stackrel{i} \to (E,E') \stackrel{p} \to B
$$
is such that    $F'$ is subspace  of $F$ and local trivializations of $E'$ 
obtained by
restricting local trivialization of $p$. If is well known (see \cite{Ha}-Theorem 
4D8 for
example) that for any ring of coefficients $\bk$ that  if  there exist cohomology
classes $c_i \in H^\ast (E,E';\bk) $ whose restrictions $\{H^\ast (i,i';\bk) 
(c_i) $
freely generate the $\bk$-module   $H^\ast (F,F';\bk)$ then the  $H^\ast
(B;\bk)$-module  $H^\ast (E,E';\bk)$ is freely generated by the $c_i$'s. In the 
case  of
fiber bundle pair $(p,p')$ with fiber $(D^k, S^{k-1})$ ($D^k$ denotes the
$k$-disk in ${\mathbb R}^k$ and its boundary $S^{k-1}$) 
an element $\tau \in H^k(E,E')$ whose restriction to $H^k(D^k, S^{k-1})$ is the 
fundamental class 
is called the {\it  Thom class of the fiber bundle pair $(p,p')$}. If the fiber 
bundle pair
$(p,p')$ with fiber $(D^k, S^{k-1})$ admits a Thom class $\tau ^{(p,p')}$ then 
the
$H^\ast(B)$-linear map 
$$
H^\ast (p)(-) \cup \tau^{(p,p')} : H^\ast(B) \to H^{\ast+k} (E,E')\,,
\quad  x \mapsto H(p)(x)\cup \tau $$
is an isomorphism, called {\it the Thom isomorphism},  and $ H^{l} (E,E' )=0
$ if
$l <k$.

 Recall, see for example \cite{Ha}-Theorem 4D10), if $\bk$ is a field of
characteristic $\neq 2$ then  every orientable fiber
bundle pair $(p,p')$ with fiber
$(D^k, S^{k-1})$ admits a Thom class. In this case  the Gysin sequence of the
sphere bundle
$p'$ can be deduced from the long exact of the pair $(E,E')$ and the Thom 
isomorphism.
In particular, one has the following commutative diagram
$$
\begin{array}{ccccccccccc}
..&H^l(E,E'  ) &\stackrel{j^{(E,E'  )}}{\to} &H ^l (E'  ) & {\to} &H ^l (E  ) 
&&\to &
H^{l+1}(E,E'  )&...\\ & \cong \uparrow&& \cong \uparrow H(p)&& || &&
&\uparrow \cong&\\ ...&H^{l-k}(B  ) &\stackrel{-\cup e^{p'}}{\to} &H ^l (B  ) &
\stackrel{H^l(p')}{\to} & H ^l (E  ) &&  \to & H^{l+1-k}( B  )&...
\end{array}
$$
where $e^{p'} = \left( H^\ast(p)\right)^{-1} \circ j^{(E,E')} 
(\tau) \in H ^{k}(B  )$ is usually called the {\it  Euler class of the sphere
bundle $p'$}.

 Let $M$ and $ N$  be  (smooth Banach and without boundary) 
 connected manifolds   and $ f : M\to N $ be a
(smooth) closed  embedding, \cite{La}-I$\!$I.$\S$2. Then we have the exact 
sequence of
fiber bundles
$$
0 \to TM \stackrel{Tf}\to TN_{|_M} \to \nu_f \to 0$$
where $TM$ and $TN$ are the tangent bundles and $\nu_f$ is  the {\it normal fiber 
bundle
of
$f$}. By
definition of an immersion, \cite{La}-I$\!$I Proposition 2.3,  this exact 
sequence split.
Hereafter we will identify
$
\nu_f
$ with a factor of
$TN_{|_M}$. 

When the fiber of $\nu_f $ is of finite dimension $k$ the embedding is called a 
{\it
$k$-embedding}. In this case the normal fiber bundle pair $(\nu _f ^D, \nu_f^S)$ 
is
defined and the  Thom class of the oriented normal fiber bundle pair
$(\pi^D,\pi^S)$,   denoted  $\tau _f$,  is called  the Thom class of the 
oriented embedding  $f$.

For any embedding $f$, the exponential map $D \subset TN \to N $ restricted  to 
$\nu_f$
is a local isomorphism on the zero cross section of the bundle $TN \to N$.  Thus 
if we
assume that $f$ is a closed embedding and that $N$ admits a partition of unity 
then, by 
\cite{La}-I$\!$V-Theorem 5.1,  there exists 

a) an open neighborhood $Z$ of the zero section of $\nu_f$,

b) an open neighborhood $U$ of $f(M)$ in $N$,

c) an isomorphism $\Theta : Z \to U$ which identifies the zero section of $\nu_f$ 
with
$f(M)$.

If moreover,  we assume that $f$ is a $k$-embedding    we 
identify $Z$ with $\nu_f ^D$  
and the isomorphism 
$$
  \nu _f ^D \stackrel{\theta}{\to} U : = \mbox{tube}\, f
$$
restricts to an isomorphism
 $$
\nu^S_f \cong \theta (\nu_f^S) := \partial \mbox{tube}\, f \,.
$$
It results from, \cite{La}-V$\!$I$\!$I-Corollary 4.2, that $\nu_f\cong \nu_f ^D$. 
The
 above discussion is summarized in the next commutative  diagram
$$
(6) \qquad \quad 
\begin{array}{ccccc}
H^\ast(N, N-f(M)) & \stackrel{j^{(N,N-f(N))}}{\longrightarrow}& H ^\ast(N)\\
\cong \downarrow \mbox{\tiny Excision} && \downarrow \\
H^\ast(\mbox{tube}\,f, \partial \mbox{tube}\,f) &
\stackrel{j^{(\mbox{\tiny tube}\,f,\partial \mbox{\tiny 
tube}\,f)}}{\longrightarrow}& H
^\ast(\mbox{tube}\,f)\\
 ||&&|| \\
H^\ast(\nu_f^D, \nu _f ^S) &
\stackrel{j^{( \nu_f^S, \nu^S _f)}}{\longrightarrow}& H
^\ast(\nu^D_f)\\
-\cup \tau _f \uparrow\cong &&\hspace{15mm} \uparrow  \pi^D \circ 
f^{-1}_{|_{f(M)}} \\
H^{\ast -k} (M) & \stackrel{ - \cap e ^{\pi^S}}{\longrightarrow}&  H^\ast(M)
\end{array}
$$
The composition
$$
 TP^\ast :  H^\ast(\nu^D_f,
\nu ^S _f) = H^\ast(\mbox{tube }f, \partial \mbox{tube} 
f)\stackrel{Excision}\cong
H^\ast(N,N-f(M)) \stackrel{j^{N, N-f(M)}}{\rightarrow } H^\ast( N)  
 $$
is called {\it the cohomology Thom-Pontryagin collapse map}.

The $H^\ast(N)$-linear
map 
$$
f^! :H ^\ast (M) \stackrel{H^\ast(\pi ^D\circ f)(-) \cup \tau _f}{\to} 
H^{\ast + k}(\nu ^D_f , \nu^S_f)
\stackrel{TP^\ast}{\to} H^{\ast+k} (N)\,.
$$
is called the {\it cohomology shriek map} and  $ f^!(1) =e_f \in H^{k}(N)$ is 
called the {\it
Euler class } of the embedding $f$.

\vspace{3mm}
\noindent{\bf Proposition 1.} {\it Let us consider the commutative diagram
$$
\begin{array}{ccc}
M &\stackrel{f}\to &N\\
\varphi \downarrow && \downarrow \psi\\
M' &\stackrel{g}\to &N'
\end{array}
$$
where 
$$
\left\{
\begin{array}{ll}
M\,,N\,, M'\,, N' \mbox{ are smooth  Hilbert manifolds (without boundary)}\\
f \,, g\,, \mbox{ are two finite codimensionnal  closed oriented embeddings}\\
 \varphi \,, \psi\ \mbox{ are two smooth maps  which preserve  orientations}\\
\hspace{10mm} \mbox{ of the
associated normal pairs} 
\end{array}
\right.
\,.
$$
Then we have the following commutative diagram
$$
\begin{array}{ccc}
H^\ast(M)&\stackrel{f^!}{\to}& H^{\ast+k}(N)\\
H^\ast(\varphi) \uparrow && \uparrow H^\ast(\psi) \\
H^\ast(M')&\stackrel{g^!}{\to} &H^{\ast+k}(N')
\end{array}
$$ }

\vspace{3mm}

\noindent {\bf  Proof.}   The maps $ \varphi$ and $\psi$ preserve orientations  
of   the normal
pairs  implies that   they have the same rank and that    the Euler classes  of 
the
embeddings are preserved.     By,
\cite{La}-I$\!$V-Theorem 6.2, there exists a vector bundle isomorphism
$\lambda
$ between
$\nu_f^D$ and the pullback $\varphi^\ast \nu_g ^D$ of $\nu_g ^D$ via $\varphi$ 
and an 
isotopy
$ h : 
\nu_f^D
\times {\mathbb R} \to \varphi^\ast \nu_g ^D$ such that $h(-,0)= f$ and 
$h(-,1)=g\circ
\lambda$. Therefore the pair $(\psi, \varphi)$ induces a well defined 
homomorphism
$$
H^\ast(\mbox{tube}\, g, \partial \mbox{tube}\, g)  \to H^\ast(\mbox{tube}\, f,
\partial
\mbox{tube}\,f ) 
$$
The remaining of the proof is now straightforward.

\hfill{$\square$}

 We have the very useful result

\vspace{3mm}

\noindent {\bf  Proposition 2.} {\it  Let $M$ and $ N$  be  (smooth Hilbert  and 
without
boundary) 
 connected manifolds   and  $f  : M \to N $ be a closed  oriented  $k$-embedding 
with Euler class
$  e_f \in H^k(N) $. The composite 
$$
H^\ast(N) \stackrel{H^\ast(f)}\to H^\ast (M) \stackrel{f^!} \to H^{\ast+k} (N) 
$$
is the right multiplication  by $e_f$. Moreover, the restriction of
$f^!$  to the image of $H^\ast(f)$ is well defined by the formula:
$$
f^! \circ H^\ast(f)=-\cup e_f\,.
$$
 In particular,    if $H_\ast(f)$ is an isomorphism then
$f^!= \left(H^\ast(f)\right)^{-1}$.}

\vspace{3mm}
\noindent{\bf Proof.} Observe that, by construction,   $ f^!$ is a $H^\ast
(N)$-module homomorphism  (here $H^\ast (M)$ is a $H^\ast(N)$-module via
$H^\ast(f)$). Then for any $x \in  H^{l}  (N)$ and any $y \in  H^\ast (M)$ 
$$
f^!(H^\ast(f)(x)\cup  y)= x\cup f^!(y)\,.
$$
 If   $\alpha  =H^\ast(f)(\beta)$, for some $\beta
\in  H^{\ast}(N)$ and if   $1$ denotes the unit in $ H^0(M)$ then 
$$f^!(\alpha )= f^!(H^\ast(f)(\beta)\cup 1 )=\beta\cup 
f^!(1)=\beta \cup e_f\,.$$ 

Since $f^!$ is well defined the right-hand part of the above formula is 
independent of
the choice of
$\beta$. This can be also established in the following way. From diagram (6) one 
deduces  
the commutative diagram
$$
\begin{array}{cccccc}
H^\ast(\nu ^D _f ,\nu ^S_f)& \stackrel{\cong}{\leftarrow}& H^\ast(N, N-f(M))&
\stackrel{j^{N,N-f(M)}}{\rightarrow}& H^\ast (N) \\
 
-\cup\tau _f \uparrow && \hspace{12mm}\uparrow -\cup \tau ^{N,N-f(M)}_f
&&
\\

H^{\ast-k}(\nu ^D _f ) = H^{\ast-k} ( M) &
\stackrel{H^\ast(f)}{\longleftarrow}&\hspace{-5mm} H^{\ast-k}(N) 
\end{array}
$$
where $\tau ^{N,N-f(M)}_f$ denotes the image of the Thom class $\tau_f$ via the
isomorphism $H^\ast(\nu ^D _f ,\nu ^S_f)\cong  H^\ast(N, N-f(M))$. Now let
$\beta\,, \beta' \in H ^\ast (N)$ be such that $H^\ast(f)(\beta)=\alpha=
H^\ast(f)(\beta')$. Thus $\beta-\beta' \in \mbox{ker}\,H^\ast(f)$ and then $ 
(\beta
-\beta') \cup \tau ^{N,N-f(M)}_f=0$. Therefore  $ (\beta
-\beta') \cup e_f=0$.

 \hfill$\square$

Let $M$ be a    simply connected compact and oriented  $m$-manifold. A slight 
modification in the
proof of  Lemma 5.4.1 in \cite{Bry} shows that the smooth map $p_0 : LM \to M$ is 
a fiber bundle.
Thus each fibration apearing in the  two
pullback diagrams 
$$
 (3) \quad \qquad 
\begin{array}{ccccccl}
 LM ^{\times 2} &\stackrel{\delta_{out}  }{\longleftarrow}&  LM \times_M 
LM&\stackrel{
\delta_{in} }{\longrightarrow}&  LM \\
\hspace{- 9mm} p_0\times p_0 \downarrow && \hspace{- 4mm}p_0 \downarrow &&
\hspace{2mm}\downarrow   (p_0, p_{\frac{1}{2}})
\\
M^{\times 2}  &\stackrel{\Delta  }{\leftarrow}&M
&\stackrel{\Delta}{\to}& M^{ \times 2}   
\end{array}
$$
are  smooth fiber bundles.

\vspace{3mm}
\noindent{\bf Theorem 3}   (\cite{C-J}-$\S 1$ and  \cite{C-G}-$\S 2$) {\it  
$\delta_{in}$ and $\delta_{out}$ are closed $m$-embeddings.}

\vspace{2mm}

  It results from Proposition 1 that we have
the  following commutative diagram
$$
 (7) \quad \qquad\qquad 
\begin{array}{ccccccl}
H^{\ast+k} ( LM     ) ^{\otimes 2} 
&\stackrel{\delta_{out}^!}\longleftarrow&
H^{\ast } ( LM \times_M LM    ) &
\stackrel{ \delta_{in}^! } \longrightarrow & H^ {\ast+k}( LM   )  \\
\hspace{- 19mm} H^\ast(p_0)^{\otimes 2} \uparrow && \hspace{4mm}H^\ast(p_0)
\uparrow &&
\hspace{2mm}\uparrow   H^\ast (p_0, p_{\frac{1}{2}})
\\
H^ {\ast+k} ( M   ) ^{\otimes 2}  &\stackrel{ \Delta^! }{\longleftarrow 
}&H^{\ast} (M  )
&\stackrel{ \Delta^!}{\longrightarrow }& H^{\ast+k}(M   ) ^{
\otimes 2}   
\end{array}
$$

In particular
$$ 
 e_{\delta _{in}} =H^\ast(p_0,p_{\frac{1}{2}})
(e_\Delta)\in H^k(LM ) \mbox{ and }  e_{\delta _{out}}= H^\ast(p_0,p_1)
(e_\Delta)
 \in H^k(LM\times LM  )
\,.
$$
where $e_\Delta$ denotes the Euler class of the diagonal embedding of $M$.

Moreover, by Proposition 1,  diagram (7) is natural with respect to orientation
preserving smooth  maps between $m$-dimensional simply connected compact and 
oriented 
manifolds.

\vspace{5mm}

 \centerline{\sc $\mathbf{\S 3}$ -  Sullivan representative of $\delta_{in}$ and 
of
$\delta_{out}$.}

\vspace{3mm}
In this section  $\bk$ is  a field of characteristic zero.
  We refer systematically the reader to \cite{FHT} for
notation,  terminology and needed results.  Nonetheless, we  precise here some 
notation
and recall some  relevant results.  If $V =\{V^i\}_{i \in
\geq 0} 
$ is a
 graded $\bk$-vector space   
$\overline{V}$ denotes the suspension of $V$,  $ (\overline{V})^n =V^{n+1}$. We 
denote
by
$
\bigwedge V$
 the free graded commutative algebra generated by $V$. A {\it Sullivan algebra} 
is a
cochain algebra of the form $ (\bigwedge V\,,d)$. It is  minimal if  $d(V) 
\subset
\bigwedge^{\geq 2}(V)$.

 Any path connected space $X$ admits a Sullivan  model,
$$
\mm _X:= (\bigwedge V,d) \stackrel{\simeq}{\to} A_{PL}(X)
$$
where $A_{PL}$ denotes the contravariant functor of piece linear differential 
forms,
\cite{FHT}-$\S 10$ and $\S 12$. If $X$ is 1-connected 
$\mm _X$ can be chosen minimal and in this case it is unique up to isomorphism.

Let $X$ and $Y$ be two path-connected spaces. Any continuous map $f : X \to Y$ 
admits a
{\it Sullivan representative}
$\mm _f :
\mm _Y
\to
\mm _X $ : there exists a diagram 
$$
\begin{array}{ccc}
A_{PL}(Y)&\stackrel{A_{PL}(f)}{\longrightarrow} &A_{PL}(X) \\
\rho _Y \uparrow \simeq && \simeq \uparrow \rho_X \\
\mm _Y & \stackrel{\mm _f}{\longrightarrow}& \mm_X
\end{array}
$$
such that $\rho_X \circ \mm_f $ and $A_{PL}(f) \circ \rho_Y $ are homotopic
homomorphisms of commutative differential graded algebras. Moreover, if $f, g : X 
\to Y$
are homotopic then $\mm _f , \mm _g : \mm_Y \to \mm_X$ are homotopic. 

Here after we will make the following identifications:
$$
H^\ast (X) = H( A_{PL}(X))=H(\mm_X)\,,\,  H^\ast (f) = H( A_{PL}(f))=H(\mm_f)\,.
$$

 It is convenient   to consider ``cofibrant'' Sullivan representative of $f$ 
called
{\it relative Sullivan model} of $f$. The map 
$f: X
\to Y$   admits also  a   relative Sullivan model,  \cite{FHT}-$\S 14$:
$$
\begin{array}{ccc}
A_{PL}(Y)&\stackrel{A_{PL}(f)}{\longrightarrow} &A_{PL}(X) \\
\rho _Y \uparrow \simeq && \simeq \uparrow \rho_X \\
\mm _Y & \stackrel{\lambda _f}{\longrightarrow}& \mm_X=(\mm_Y \otimes \bigwedge 
W,d)
\end{array}
$$
where $\lambda_f$ denotes the natural inclusion of $\mm _Y $ into $ \mm_Y \otimes 
\bigwedge
W$ and $d$ is such that $de_i \in \mm_Y \otimes \bigwedge (e_1,e_2,...e_{i-1})$ 
for some
linear basis $e_1,e_2,..., e_n,....$ of $W$.

It important to observe here that :

1)  the relative Sullivan model of $f$ is unique up to isomorphism 
\cite{FHT}-Theorem 
14.12. 

2)   all the rational homotopy informations on $f$ are hidden in the differential 
$d$.

\vspace{3mm}

\noindent{\bf  Theorem 4.}{ \it Let $M$ be a 1-connected compact oriented 
manifold and  
let $\mm _M = (\bigwedge V,d)$ be  a minimal model of
$M$ whose product is denoted by $\mu : \bigwedge V \otimes \bigwedge V \to 
\bigwedge V$.

a) A Sullivan representative of $\delta_{out}$ is explicitly given  by:
$$
\begin{array}{ccc}

(\bigwedge V \otimes \bigwedge \bar V, d)^{\otimes 2}\cong   ((\bigwedge 
V)^{\otimes
2} 
\otimes 
(\bigwedge \bar V)^{\otimes 2} , d)& 
\stackrel{ \mu \otimes id\otimes id }{\longrightarrow} & (\bigwedge
V
\otimes
(\bigwedge
\bar V)^{
\otimes 2}, d)\\
||&&||
\\
\mm _{LM} ^{\otimes 2} \cong \mm _{LM\times LM} & \stackrel{\mm _{\delta_{out}}
}{\longrightarrow} &\mm _{LM\times_M LM}
\end{array}
\,.
$$

b)  A Sullivan representative of $\delta_{in}$ is explicitly given  by:
$$
\begin{array}{ccc}
    ((\bigwedge V)^{\otimes
2} 
\otimes 
(\bigwedge \bar V)^{\otimes 2} , d') 
&\stackrel{  \mu\otimes id\otimes id  } {\longrightarrow}
&(\bigwedge
V
\otimes
(\bigwedge
\bar V)^{
\otimes 2}, d)\\
|| &&|| \\
 \mm_{LM}' 
&\stackrel{\mm _{\delta _{in}}  } {\longrightarrow}
&\mm _{LM\times_M LM}
\end{array}
$$
 The differentials
appearing
 in the above  statement  will be  precisely described during the proof.}

Observe here that while  the homomorphisms of graded algebras $\mm
_{\delta _{in}}$ and $ \mm _{\delta _{out}}
$ are equal the homomorphism $ H(\mm _{\delta _{in}} ) $ and $ H(\mm _{\delta 
_{out}} )$
are different.   The remaining of this section is devoted to the proof of Theorem 
4.
The proof follows from a sequence of lemmas which rely  on the ``preservation
under pushout'' property  of the relative Sullivan models:

\vspace{3mm}

\noindent{\bf  Lemma 5.} \cite{FHT}-Proposition 15.8  { \it  A  pullback diagram 
in
the category of topological spaces where $p$ is a fibration :
$$
\begin{array}{cccc}
f^\ast E &\stackrel{f'}{\to}&E\\
p'\downarrow && \downarrow  p\\
B' &\stackrel{f} {\to} & B
\end{array}
$$ 
converts  into the  pushout diagram in the category of commutative differential 
graded
algebras
$$
\begin{array}{cccc}
\mm_B & \stackrel{\mm_f}{\to} &\mm_{B'}\\
\lambda_p \downarrow && \downarrow \lambda _{p'}\\
\mm_E =(\mm_B \otimes \bigwedge W, d)& \stackrel{\mm_{f}\otimes id}{\to} 
&(\mm_{B'}
\otimes
\bigwedge W, d')=\mm_{E'}
\end{array}
$$
}

Let us just precise here that the differential  $d'=\varphi \circ d\circ \varphi 
^{-1} $
 where  $\varphi$ is 
 the canonical isomorphism of commutative graded algebras  such that the next 
diagram
commutes
$$
\begin{array}{cccc}
\mm_{B'} \otimes _{\mm _B} \left( \mm_B \otimes \bigwedge W\right)& \stackrel
{\varphi}{\longrightarrow} &\mm_{B'} \otimes \bigwedge W & \qquad \quad a\otimes 
b
\otimes c
\mapsto a\mm_f(b)\otimes c\\
\mm_f\otimes id\otimes id  \uparrow &&\uparrow \mm_f \otimes id\\
\left( \mm_{B} \otimes _{\mm _B}  \mm_B \right)  \otimes \bigwedge W& \stackrel
{\mu \otimes id }{\longrightarrow} &\mm_{B} \otimes \bigwedge W & \qquad \quad 
a\otimes b
\otimes c
\mapsto ab\otimes c
\end{array}
\,.$$

\vspace{3mm}

\noindent{\bf  Lemma 6.} \cite{FHT}-$\S 15$-Example 1 { \it The relative Sullivan 
model 
of the   fibration $ p_0\times p_1  : M ^I \to M\times M $ is of the form  
$$
\mm_M^{\otimes 2}= (\bigwedge V, d)^{ \otimes 2} 
\stackrel{\lambda_{p_0\times p_1} }{\hookrightarrow} (
(\bigwedge V)^{\otimes 2} 
\otimes \bigwedge \bar V, d)=\mm_{M^I }\,.
$$
}

Here the differential is precisely defined as follows.\\
$
\left\{ 
\begin{array}{ll}  d (v\otimes 1\otimes 1) &=dv\otimes 1\otimes \bar 1 \\
 d (1\otimes v\otimes \bar 1) &=1\otimes dv\otimes \bar 1 \\
  d (1\otimes 1\otimes \bar v)& = (v\otimes 1-1\otimes v)\otimes \bar 1 -
\sum_{i=1}^\infty
\frac{(sd)^i}{i!}(v\otimes 1\otimes \bar 1)
\\
\end{array} \right.\,, \qquad  v\in V\,,  \bar v \in \bar V
$  and $ s$ being the
unique degree -1 derivation   of $  \bigwedge V \otimes \bigwedge V
\otimes \bigwedge \bar V$ defined by 

$\left\{
\begin{array}{lll}
s(v\otimes 1\otimes \bar 1)=1\otimes 1\otimes \bar v =
s(1\otimes v\otimes \bar 1)\\
s(1\otimes 1\otimes \bar v)=0
\end{array} \right. \,, \qquad v \in V \,, \bar v \in \bar V\,.
$

 We denote by 1 (respectively  $\bar 1$) the unit of $\bigwedge V$ (respectively  
$\bigwedge \bar V$).

\vspace{3mm}
\noindent{\bf Lemma 7.} {\it  The relative Sullivan model  for $p_0:   LM \to M $ 
is 
$ \lambda_{p_0}  :
\mm_M= (\bigwedge V,d) \hookrightarrow  (\bigwedge V\otimes \bigwedge \bar V, 
 d)=\mm_{LM}\,.
$}

\vspace{3mm}
\noindent{\bf Proof } Since  multiplication $\mu : (\bigwedge 
V,d)^{ \otimes 2}  \to (\bigwedge V,d)$ is a Sullivan representative of
$\Delta$, by Lemma 6 and Lemma 7   the pullback diagram 
$$
\begin{array}{ccl}
  LM
&\stackrel{}{\to}& M^I\\
p_0 \downarrow &&
 \hspace{2mm}\downarrow  p_0\times p_1
\\
M
&\stackrel{\Delta}{\to}& M^{ \times 2}   
\end{array}
$$
yields the pushout diagram in the category of
commutative differential graded algebras
$$
\begin{array}{ccl}
  (\bigwedge  V,d) ^{\otimes 2}  &\stackrel{\mu} \to &(\bigwedge V,d)\\

\lambda_{p_0\times p_1}  \downarrow &&
 \hspace{5mm}\downarrow  \lambda_{p_0}  
\\
((\bigwedge V)^{\otimes 2}  
\otimes \bigwedge \bar V, d)
& \stackrel{\mu \otimes id} \to& (\bigwedge V \otimes \bigwedge \bar V, \bar d)  
\end{array}\,.
$$
The relative Sullivan model  for $p_0:   LM \to M $ is $\lambda_{p_0}$ and 
$
\left\{ 
\begin{array}{ll}
 d (1\otimes v)&=dv\otimes \bar 1 \\ 
 d(1\otimes  \bar v)& = -s(dv\otimes \bar 1) \\
\end{array} \right.
$  with $  s$ being the
unique degree -1 derivation of $\bigwedge V \otimes \bigwedge \bar V$ defined by 
$ 
\left\{ \begin{array}{ll} s(v\otimes \bar 1)&=1 \otimes \bar v \\ s(1\otimes \bar
v)&=0\end{array}
\right.$, $
 v
\in V\,, \bar v \in \bar V\,.$

\hfill $\square$

\vspace{3mm}
\noindent{\bf Proof of part a) of theorem 4.} By Lemma 5 and Lemma 7 the pullback 
diagram of fibrations
$$
\begin{array}{ccccccl}
  LM \times _M LM &\stackrel{
 \delta_{out} }{\to}&  LM \times LM\\
\hspace{-1mm}p_0 \downarrow &&
\hspace{9mm}  \downarrow p_0\times  p_0\\
M
&\stackrel{\Delta}{\to}& M^{ \times 2}  
\end{array}
$$
converts into the pushout
diagram 
$$
\begin{array}{ccc}
 (\bigwedge V,d)^{\otimes 2} & \stackrel{\mu}{\to} &(\bigwedge V,d) \\
\lambda_{p_0\times p_0}  \downarrow && \downarrow  \lambda_{p_0}  \\
 \mm
_{LM\times LM}  = ((\bigwedge V)^{ \otimes 2}   \otimes (\bigwedge \bar 
V)^{\otimes 2} ,
d)  &\stackrel{\mu\otimes id \otimes id }{\to}& (\bigwedge V \otimes
(\bigwedge
\bar V)^{ 
\otimes 2} , d)=\mm
_{LM\times_M LM}
\end{array}
$$
where 
$$
\left\{ \begin{array}{ll}
d (v \otimes \bar 1 \otimes \bar 1) &= dv \otimes \bar 1 \otimes \bar 1 \\
d(1 \otimes \bar v \otimes \bar 1)& =-s(dv  \otimes \bar 1 \otimes \bar 1) \\
d(1 \otimes \bar 1 \otimes \bar v)& =-s'(dv  \otimes \bar 1 \otimes \bar 1)   
 \end{array} \right. \,, v \in V \,, \overline{v} \in \overline{V}\,.
$$
Here  $ s$ and $s'$ are  the
unique degree -1 derivation of $\bigwedge V 
\otimes (\bigwedge \bar V)^{\otimes 2}  $ defined by 
$$
\left\{ \begin{array}{ll}
s (v  \otimes \bar 1 \otimes \bar 1) &= 1 \otimes \bar v \otimes \bar 1 \\
s' (v  \otimes \bar 1 \otimes \bar 1) &= 1 \otimes \bar 1 \otimes \bar v \\
s(1 \otimes \bar v \otimes \bar 1)& =0= s'(1 \otimes \bar v \otimes \bar 1)\\
s(1 \otimes \bar 1 \otimes \bar v)& =0 =s'(1 \otimes \bar 1 \otimes \bar v)
 \end{array} \right.
 \,, \qquad  v\in V\,, 
 \bar v \in \bar V\,.
$$

Thus  $
\mm _{\delta _{out}} = \mu\otimes id\otimes id $.

\hfill $\square$

\vspace{3mm}
\noindent{\bf Lemma 8.} {\it A relative Sullivan model for $ (p_{0},p_ 
{\frac{1}{2}}): LM
\to M\times M  $ is 
$$
  \lambda_{p_0\times p_\frac{1}{2}} : (\bigwedge
V,d)^{\otimes 2}  \hookrightarrow (( \bigwedge V)^{ \otimes 2}  
\otimes (\bigwedge \bar V) ^{\otimes 2} , d')=\mm _{LM}'\,. 
$$}

\noindent{\bf Proof.} By Lemma 6 and  Lemma 7  
the  pullback diagram 
$$
 (4) \qquad \qquad 
\begin{array}{lcl}
 LM &\stackrel{\delta' }{\to}& {M^{[0,\frac{1}{2}]}}\times {M^{[\frac{1}{2}, 
1]}}\\
\hspace{-9mm} (p_0 , p_{\frac{1}{2}})  \downarrow   &&\hspace{1mm} \downarrow
p_0\times p_\frac{1}{2}\times p_\frac{1}{2}\times p_1
\\
 M^{\times 2}  &\stackrel{\Delta'  }{\rightarrow}& M^{\times 4}
\end{array}\,, 
$$ 
  converts into the pushout diagram 
$$
\begin{array}{ccc}
(\bigwedge V,d)^{\otimes 4} 
&\stackrel{\mm_{\Delta'}}{\to}& (\bigwedge V,d)^{\otimes 2} \\
\lambda _{p_0\times p_{\frac{1}{2}}}^{ \otimes 2}  
\downarrow && \downarrow \lambda _{p_0\times p_{\frac{1}{2}}} \\ 
(\bigwedge V \otimes \bigwedge
V \otimes \bigwedge \bar V, d)^{\otimes 2}   & \stackrel{\mm _{\delta '}}{\to} &(
(\bigwedge V )^{\otimes 2}  
\otimes (\bigwedge \bar V)^{\otimes 2} , d')
\end{array}
$$
where $\mm_{\Delta '} (a\otimes b\otimes a'\otimes b')= (-1)^{|b'|(|a'|+|b|)} 
ab'\otimes
ba'$, with $\lambda_{p_0\times p_{\frac{1}{2}}}$ defined by\\ $
\left\{
\begin{array}{ll} v\otimes 1 \mapsto v \otimes 1\otimes \bar 1 \otimes \bar 1\\
1\otimes v \mapsto 1 \otimes v \otimes \bar 1 \otimes \bar 1 
\end{array}
\right.
\,, 
$ and $d'$ defined by 
$$
\left\{ 
\begin{array}{ll} 
d'(v\otimes 1\otimes \bar 1 \otimes \bar 1)&=dv \otimes 1 \otimes \bar 1 \otimes 
\bar 1
\\ d'(1\otimes v\otimes \bar 1 \otimes \bar 1)&=1 \otimes dv \otimes \bar 1 
\otimes \bar
1 
\\ 
d'(1\otimes 1\otimes \bar v \otimes \bar 1)&=(v \otimes 1 - 1\otimes v) \otimes 
\bar 1
\otimes \bar 1 - \sum _{n \geq 1} \frac{(sd')^n}{n!}(v\otimes 1\otimes \bar 1 
\otimes
\bar 1)\\
d'(1\otimes 1\otimes \bar 1\otimes \bar v)&=(1 \otimes v - v\otimes 1) \otimes 
\bar 1
\otimes \bar 1 - \sum _{n \geq 1} \frac{(sd')^n}{n!}(1\otimes v\otimes \bar 1 
\otimes
\bar 1)\\

\end{array} \right.
$$
where  $ s'$ denotes the 
unique degree -1 derivation of $(\bigwedge V)^{\otimes 2}  
\otimes (\bigwedge \bar V)^{\otimes 2}  $ defined by 
$$
\left\{ 
\begin{array}{ll} 
s'(v\otimes 1\otimes \bar 1 \otimes \bar 1)&=1 \otimes 1 \otimes \bar v \otimes 
\bar 1 \\
s'(1\otimes v\otimes \bar 1 \otimes \bar 1)&=1 \otimes 1 \otimes \bar 1 \otimes 
\bar v 
\\ 
s'(1\otimes 1\otimes \bar v \otimes \bar 1)&=0\\
s'(1\otimes 1\otimes \bar 1\otimes \bar v)&=0\\
\end{array} \right.
 \,, \qquad  v\in V \,,  \bar v \in \bar V 
\,.$$
\hfill $\square$

\vspace{3mm}
\noindent{\bf Lemma 9.} {\it  Let $\mu$ (respectively  $\bar \mu$) be 
the product of $\bigwedge V$ (respectively  of $\bigwedge \bar V $). The 
surjective
homomorphism of differential graded algebras 
$$
\mu \otimes \bar \mu  :\mm _{LM}'= ((\bigwedge V)^{\otimes 2}  \otimes
(\bigwedge
\bar V)^{ \otimes 2} , d')\to
(\bigwedge V\otimes \bigwedge \bar V,d) =\mm_{LM},.
$$ 
is a quasi-isomorphism. }

\vspace{3mm}
\noindent{\bf Proof.}  Formulae $
 \left\{
\begin{array}{lll} 
v\otimes 1 \otimes \bar 1 \otimes \bar 1  \mapsto (v\otimes 1 + 1\otimes v)  
\otimes
\bar 1 \otimes \bar 1\\

1\otimes v \otimes \bar 1 \otimes \bar 1    \mapsto (v\otimes 1 - 1\otimes v) 
\otimes \bar 1 \otimes \bar 1\\

1\otimes 1 \otimes \bar v \otimes \bar 1  \mapsto 1\otimes 1\otimes (\bar
v\otimes 1 + 1\otimes \bar v)  \\

1\otimes 1 \otimes \bar 1 \otimes \bar v  \mapsto 1\otimes 1\otimes (\bar
v\otimes 1 - 1\otimes \bar v)  \\

 \end{array}
\right.
\,,\qquad v \in V\,, \bar v \in \bar V
$
define an automorphism of commutative graded algebras of  $(\bigwedge V)^{\otimes 
2} 
 \otimes ( \bigwedge \bar V)^{ \otimes 2} $. This automorphism
exhibits the differential graded algebra $((\bigwedge V)^{\otimes 2} 
 \otimes ( \bigwedge \bar V)^{ \otimes 2} ,d') $ as the tensor product
 of $\mm _{LM}$ with
an acyclic commutative differential graded algebra.  Therefore the surjective
homomorphism of differential graded algebras 
$$
\rho :\mm _{LM}'= ((\bigwedge V)^{\otimes 2} 
 \otimes ( \bigwedge \bar V)^{ \otimes 2} ,d')\to
(\bigwedge V\otimes \bigwedge \bar V,d) =\mm_{LM},.
$$ 
uniquely defined by formulae 
$
 \left\{
\begin{array}{lll} 
v\otimes 1 \otimes \bar 1 \otimes \bar 1  \mapsto  v   \otimes \bar 1\\

1\otimes v \otimes \bar 1 \otimes \bar 1  \mapsto   v\otimes   \bar 1\\

1\otimes 1 \otimes \bar v \otimes \bar 1  \mapsto 1\otimes   \bar v \\

1\otimes 1 \otimes \bar 1 \otimes \bar v  \mapsto 1\otimes   \bar v \\

 \end{array}
\right.
\,,\qquad v \in V\,, \bar v \in \bar V
$ 
is a quasi-isomorphism.  By construction $\rho = \mu \otimes \bar \mu$. 

\hfill $\square$

\vspace{3mm}
\noindent{\bf Proof of part b) of theorem 4.}  By Lemma 5  and Lemma 8 the
pullback  diagram of fibrations
$$
\begin{array}{ccccccl}
  LM \times _M LM &\stackrel{
 \delta_{in} }{\to}&  LM\\
\hspace{-1mm}p_0 \downarrow &&
\hspace{9mm}  \downarrow p_0\times  p_\frac{1}{2}\\
M
&\stackrel{\Delta}{\to}& M^{ \times 2}  
\end{array}
$$
converts into the pushout
diagram 
$$
\begin{array}{ccc}
 (\bigwedge V,d)^{\otimes 2} & \stackrel{\mu}{\to} &(\bigwedge V,d) \\
\lambda_{p_0\times p_\frac{1}{2}}  \downarrow && \downarrow  \lambda_{p_0}  \\
 \mm_{LM}'= ((\bigwedge V)^{ \otimes 2}  \otimes (\bigwedge \bar V)^{\otimes 2} , 
d') 
&\stackrel{\mu \otimes id\otimes id  }{\to}& (\bigwedge V \otimes
(\bigwedge
\bar V)^{
\otimes 2} , d)=\mm
_{LM\times_M LM}
\end{array}\,.
$$
Thus  $
\mm _{\delta _{in}} = \mu\otimes id\otimes id $.

\hfill $\square$

\vspace{5mm}

 \centerline{\sc $\mathbf{\S 4}$ -  Proof of theorem A.}

\vspace{3mm}
In this section $\bk$ is a  field of characteristic zero and we use notation 
introduced in
section $\S$4.

First recall, \cite{Br}-V$\!I$-Theorem 12.4,  that the
Euler  class of the diagonal embedding  $\Delta : M \to M
\times M
$ (also called the diagonal class) is the cohomology class
$$ e_\Delta = \sum _{l} (-1) ^{|\beta _j|}  \hat \beta _l \otimes \beta
_l
\in H^m (M\times M ) = \left( H^\ast (M) \otimes H^\ast(M)\right) ^m
$$
where $\{\beta _l \}$ denotes a homogeneous linear basis of $H^\ast (M)$ and 
$\{\hat
\beta _l\}$ its Poincar\'e dual basis.

Let  $b_l$ (respectively  $\hat b_l$) be a cocycle of $(\bigwedge
V,d) $ representing 
$\beta _l \in H^ \ast(\bigwedge V,d)= H ^\ast(M)$  
(respectively  $\hat \beta _l \in H ^\ast(\bigwedge
V,d)= H^\ast(M)$).  Since, 
$
e_\Delta =[\sum _l (-1)^{|b_l|}  \hat b_l \otimes b_l]\in
H^m (M\times M)
$
it results from  Theorem 3  and the definition of  
$\lambda_{p_0,p_{\frac{1}{2}}}$
given in  lemma 9,  that 
 $$
e_{\delta_{in}} = [ \lambda_{p_0\times p_\frac{1}{2}} (\sum _l
(-1)^{|b_l|}\hat b_l
\otimes b_l)]=[\sum _l (-1)^{|b_l|} 
\hat b_l \otimes b_l\otimes \bar 1 \otimes \bar 1] \in H^m{ LM}\,.
$$

We consider the  {\bf degree $\bold m$ linear map} 
$$
\begin{array}{r}
 \mm _{LM\times _MLM} = (\bigwedge V \otimes (\bigwedge \bar V)^{\otimes 2} ,d)
\stackrel{ \mm_{\delta_{in}}^! }\longrightarrow 
((\bigwedge V)^{\otimes 2} \otimes  (\bigwedge V)^{ \otimes2} ,d') =\mm_{LM}' \\
 (a \otimes \bar x \otimes \bar y)
\mapsto  
\hspace{6mm}\left( \frac{1}{2}(a\otimes 1 +
1\otimes a)\otimes \bar x \otimes \bar y \right) \sum _l (-1)^{|b_l|} 
\hat b_l \otimes b_l\otimes \bar 1 \otimes \bar 1
\end{array}\,.
$$
 It is 
straightforward to check that  $d ' \circ \mm _{\delta_{in}}^!= (-1)^{m} \mm
_{\delta_{in}}^!\circ d$.

 Now  we can state precisely  the first part of theorem A.

\vspace{3mm} 

\noindent{\bf Theorem  10.} {\it  With notation  introduced above, the 
composition
$$
\mm _{LM} ^{\otimes 2}\cong \mm _{LM\times LM } \stackrel{\mm_{\delta _{in}}^! 
\circ
\mm_{\delta _{out}}} \longrightarrow  \mm ' _{LM} \stackrel{\mu\otimes \bar 
\mu}\to \mm
_{LM}\,,
$$
 where $\mu\otimes \bar \mu $ is considered in Lemma 10,  induces the  dual of 
the  loop
coproduct 
$
\Phi ^\vee : H^\ast (LM) \otimes H^\ast(LM)
\to H^{\ast +m}(LM)
$.}

\vspace{3mm}
\noindent{\bf Proof.}    Let  $ \phi :  
\mm _{LM} ^{\otimes 2}\cong \mm _{LM\times LM }   \longrightarrow  \mm ' _{LM}
$ be  a degree $m$ $\bk$-linear map such that $H(\phi)= \delta_{in}^!$.  Then by,
Proposition 2,
$H(\phi\circ  \mm_{\delta_{in }}) = - \cup e_{\delta_{in}}$.  On the otherhand, 
since $
H(\phi
\circ
\mm_{\delta_{out}})= \Phi^\vee$, in order to prove theorem 10 it suffices to 
prove that 
$$
\mbox{im}(H^\ast(\delta_{out})) \subset  \mbox{im}(H^\ast(\delta_{in})) 
\,.
$$
To prove this inclusion, we consider the commutative diagram
$$
\begin{array}{ccc}
\mm_M^{\otimes 2} &\stackrel{\lambda_{0,\frac{1}{2}}} \rightarrow &\mm_{LM}'\\
\lambda_{0} ^{\otimes 2}  \downarrow && \simeq \downarrow \mu \otimes \bar \mu \\
\mm_{LM}^{\otimes 2} &\stackrel{\hat \mu}\rightarrow &\mm_{LM}
\end{array}
$$
where $\hat \mu$ is the product in the graded algebra $\mm_{LM}$. By Lemma 9,
$\mu\otimes \bar \mu $ is a surjective quasi-isomorphism, thus Lemma 14.4 in 
\cite{FHT}
implies that there exists $\theta : \mm_{LM}^{\otimes 2} \to \mm'_{LM} $ such 
that 
$\theta \circ \lambda_{0}^{\otimes 2} = \lambda_{0,\frac{1}{2}}$ and $\mu \otimes 
\bar
\mu \circ \theta = \hat \mu$. Therefore if $\alpha $ is a cocycle in 
$\mm_{LM}^{\otimes
2}$ then $\theta (\alpha)$ is a cocyle of $ \mm_{LM} '$ and 
$\mm_{\delta_{out}}(\alpha)
= \mm_{\delta_{in}}(\theta(\alpha))$.

\hfill{$\square $}

\noindent{\bf End of the proof of   Theorem A. }

The naturality assumptions in Theorem A for the loop coproduct  follows from 
those in
Propositions 1 and the fact that two Sullivan representatives $ \varphi$ and
$\psi$  of a map $ f : X
\to Y$ are  identified, at the homology level via isomorphisms :
$$ 
\begin{array}{ccc}
H(\mm _Y) &\stackrel{H(\varphi)}\to & H(\mm_X)\\
||&&||\\
H(\mm _Y) &\stackrel{H(\psi)}\to & H(\mm_X)
\end{array}
\,.
$$

The naturality for the loop product follows, in the same way,  directly from the 
construction
performed in  \cite{FTV2}-3.5.

\hfill{$\square$}

\vspace{5mm}

 \centerline{\sc $\mathbf{\S 5}$ -  Example $M= {\mathbb C}P^n$.}

\vspace{3mm}

One know that  $\mm_ {{\mathbb C}P^n}= (\bigwedge (x,y) ,d)$ with $|x|=2$, $dy=x 
^{n+1}$ and
that $H^\ast({\mathbb C}P^m)= \bigwedge x/(x^{m+1})$. Moreover there is a 
quasi-isomorphism
$  (\bigwedge (x,y) ,d) \stackrel{\simeq}{\twoheadrightarrow} (\bigwedge 
x/(x^n+1),0)= \bar 
\mm_ {{\mathbb C}P^n}$.   For
simplicity we will replace  $ \mm_ {{\mathbb C}P^n}$ by  $\bar \mm_ {{\mathbb
C}P^n}$ so that  we obtain: 
$$
\begin{array}{rll}
\bar \mm _{LM} &= (\bigwedge x/(x^n+1) \otimes\bigwedge (\bar x ,\bar y), \bar 
d)\\

\bar \mm _{LM\times _MLM } &= (\bigwedge x/(x^n+1) \otimes\bigwedge (\bar x ,\bar 
y) \otimes\bigwedge
(\bar x ,\bar y) , \bar d)\\

\bar \mm' _{LM} &= (\bigwedge x/(x^n+1) \otimes \bigwedge x/(x^n+1) 
\otimes\bigwedge (\bar x ,\bar y)
\otimes\bigwedge (\bar x ,\bar y) , {\bar d}')\\

\bar \mm _{\delta _{in}}^! (x^k \otimes \bar u \otimes \bar v)&= 
\frac{1}{2}\left(
\sum_{i=0}^n x^i+k \otimes x^{n-i} + x^i \otimes x^{n-i+k}\right)\otimes \bar u 
\otimes
\bar v\\

\bar \mm _{\delta _{out}} (x^k \otimes x^l \otimes \bar u \otimes \bar 
v)&=x^{k+l}
\otimes \bar u \otimes\bar v\\

 \rho\circ \bar \mm _{\delta _{in}}^!\circ \bar \mm _{\delta _{out}}
 (x^k \otimes x^l
\otimes \bar x^\epsilon \bar y ^s  \otimes  \bar x^{\epsilon'} \bar y ^{s'})
&= 
\left\{
\begin{array}{ll}
0 &\mbox{ if } k+l >0\\
x^n \otimes \bar x ^{\epsilon + \epsilon'}\bar y^{s+s'}  &\mbox{ if } k=l =0 \,.
\end{array}
\right.
\end{array}
$$
Where $\epsilon$ and $\epsilon '$ is   $0$ or $1$. Since, 
$$
\{
1,  [x \otimes  \bar y^s], [x^2 \otimes  \bar y^s], ..., [x^n \otimes  \bar y^s],
[x \otimes \bar x  \bar y^s], ..., [x^{n-1} \otimes \bar x  \bar y^s]\}_{s\geq 0}
$$
is a linear basis of $H^\ast (LM)$ we have  obtained an explicit formula for 
$\Phi^\vee$
in this basis.

Let us denote by $\theta : \bar \mm _{LM} \to \bar \mm _{LM}\otimes \bar \mm 
_{LM}$ a
map inducing the dual product. From the formula obtained for $\theta$ in 
\cite{FTV2}-3.6
one deduce also that 
$$
\begin{array}{ll}
\theta \circ 
 \rho\circ \bar \mm _{\delta _{in}}^!\circ \bar \mm _{\delta _{out}}
 (x^k \otimes x^l
\otimes \bar x^\epsilon \bar y ^s  \otimes  \bar x^{\epsilon'} \bar y ^{s'})\\
\qquad \qquad = 
\left\{
\begin{array}{ll}
(1\otimes \bar x + \bar x \otimes 1 )\sum _{j=0}^{s+s'}
 \left( \begin{array}{l} s+s'\\j \end{array} \right) \,
 x^n \otimes \bar y ^i \otimes x^n
\otimes \bar y ^{s+s'-j}  &\mbox{ if } \epsilon +\epsilon' = 1\\
= \sum _{j=0}^{s+s'}
 \left( \begin{array}{l} s+s'\\j \end{array} \right) \,
 x^n \otimes \bar y ^i \otimes x^n
\otimes \bar y ^{s+s'-j}  &\mbox{ if } \epsilon =\epsilon' = 0\,.
\end{array}
\right.
\end{array}
$$

This proves in particular that the composition $\Phi \circ P$ is not
trivial.

\vspace{5mm}

 \centerline{\sc $\mathbf{\S 6}$ -  Proof of theorem B.}

\vspace{3mm}

Let us denote by $\bigwedge V \otimes \bigwedge ^p \bar V$ the subvector space 
of  $\bigwedge V \otimes
\bigwedge \bar  V$ generated by the words of length $p$ in $\bar V$. The 
differential of
$ \mm _{LM}=(\bigwedge V \otimes \bigwedge \bar V, d) $ satisfies $
d ( \bigwedge V \otimes \bigwedge ^p \bar  V) \subset  \bigwedge V \otimes 
\bigwedge ^p
\bar V$. M. Vigu\'e has proved,  \cite{V}-Theorem 3.3,  that

\centerline{$
H^n_{(i)} (LM) = H ^n(  \bigwedge V \otimes \bigwedge ^i V, d)\,. $}

To achieve the proof of  theorem B it remains to prove  that the  maps
$\mm_{\delta_{in}}^!$, $\mm_{\delta_{out}}$, and $\mu \otimes \bar \mu $   
respect the
word length in $\bar V$. The obvious relation

\centerline{$
\mu \otimes \bar \mu \left( \bigwedge^r V \otimes \bigwedge^t V \otimes 
\bigwedge^p \bar 
V
\otimes
\bigwedge^q
\bar  V 
\right) 
\subset \bigwedge^{r+t}  (V ) \otimes \bigwedge^{p+q} \bar  V  
\,.
$}

proves that $\mu\otimes \bar \mu $ respects the word length in $\bar V$. For 
$\mm_{\delta_{in}}^!$,
$\mm_{\delta_{out}}$  (defined in Theorem 10 and Theorem 4) the same argument 
works. 

\hfill{$\square$}

\hspace{-1cm}\begin{minipage}{19cm} \small
david.chataur@univ-angers.fr    and jean-claude.thomas@univ-angers.fr\\
  D\'epartement de math\'ematique  \\
Facult\'e des Sciences  \\
2, Boulevard Lavoisier\\
 49045 Angers, France .
\end{minipage}

\end{document}